\documentclass[a4paper,11pt]{article}
\usepackage{cancel}
\usepackage{amsmath}
\usepackage{amssymb}
\usepackage{setspace}
\usepackage{bm}
\usepackage{cite}
\usepackage{breakcites}

\usepackage{soul}
\usepackage{rotating}
\usepackage{adjustbox}

\usepackage{arydshln}
\usepackage{adjustbox, lipsum}
\usepackage{enumerate}
\usepackage[ruled,vlined,linesnumbered]{algorithm2e}
\usepackage{rotating}
\usepackage[utf8]{inputenc}
\usepackage{amsfonts}
\usepackage{pstricks}
\usepackage{pdflscape}
\usepackage{graphicx}
\usepackage{soul}
\usepackage{afterpage}
\usepackage{graphics}

\usepackage{rotating}
\usepackage{natbib}
\usepackage{booktabs}

\usepackage{authblk}
\usepackage{multirow}
\usepackage{mathtools}
\usepackage{scalerel}
\newtheorem{prop}{\bf Proposition}[section]

\newtheorem{remark}{Remark}[section]
\SetKwRepeat{Do}{do}{while}
\providecommand{\norm}[1]{\lVert#1\rVert}
\DeclarePairedDelimiter\abs{\lvert}{\rvert}%
\newcommand{\convex}{\ensuremath{\mbox{C-OWA-SVM}}\xspace}
\newcommand{\kernelconvex}{\ensuremath{\mbox{C-OWA-SVM}_K}\xspace}
\newcommand{\nonconvex}{\ensuremath{\mbox{NC-OWA-SVM}}\xspace}
\newcommand{\nonconvexz}{\ensuremath{\mbox{NC-OWA-SVM}(\hat{\bm{z}})}\xspace}
\newcommand{\nonconvexdual}{\ensuremath{\mbox{NC-OWA-SVM}_D(\hat{\bm{z}})}\xspace}
\newcommand{\nonconvextilde}{\ensuremath{\mbox{NC-OWA-SVM}_D}\xspace}
\newcommand{\kernelnonconvex}{\ensuremath{\mbox{NC-OWA-SVM}_K}\xspace}
\newcommand{\MProblem}{\ensuremath{\mbox{UB}_{\mbox{M}}}\xspace}
\newcommand{\bUB}{\ensuremath{\mbox{UB}_{\mbox{b}}}\xspace}
\newcommand{\bLB}{\ensuremath{\mbox{LB}_{\mbox{b}}}\xspace}

\newcommand{\dem}{\par \noindent{\bf Proof:} }
\newcommand{\fin}{\hfill $\square$  \par \bigskip}

\usepackage{xcolor}
\definecolor{gr}{rgb}{0,0.5,0}
\definecolor{armygreen2}{rgb}{0.99, 0.53, 0.93}
\definecolor{armygreen}{rgb}{0.19, 0.53, 0.43}

\author[(a)]{Alfredo Marín}
\author[(b),(c)]{Luisa I. Martínez-Merino}
\author[(b),(c)]{Justo Puerto}
\author[(d)]{Antonio M. Rodríguez-Chía}

\affil[(a)]{\small{Departamento de Estadística e Investigación Operativa, Facultad de Matemáticas, Universidad de Murcia, Murcia, Spain, amarin@um.es}}
\affil[(b)]{\small{Departamento de Estadística e Investigación Operativa, Facultad de Matemáticas, Universidad de Sevilla, Sevilla, Spain}}
\affil[(c)]{\small{Instituto de Matemáticas de la Universidad de Sevilla (IMUS), Sevilla, Spain, lmmerino@us.es, puerto@us.es}}
\affil[(d)]{\small{Departamento de Estadística e Investigación Operativa, Facultad de Ciencias, Universidad de Cádiz, Puerto Real (Cádiz), Spain, antonio.rodriguezchia@uca.es}}

\title{The soft-margin Support Vector Machine with ordered weighted average}
\date{}

\allowdisplaybreaks
\setcitestyle{authoryear,open={(},close={)}}

\makeatletter
\let\oldabs\abs
\def\abs{\@ifstar{\oldabs}{\oldabs*}}
\makeatother
\begin{document}

%
\maketitle

\begin{abstract}
	This paper deals with an extension of the Support Vector Machine (SVM) for classification problems where, in addition to maximize the margin, i.e., the width of strip defined by the two supporting hyperplanes, the minimum of the ordered weighted sum of the deviations of miclassified individuals is considered. Since the ordered weighted sum includes as particular case the sum of these deviations, the classical SVM model is a particular case of the model under study.
	A quadratic continuous formulation for the case in which weights are sorted in non-decreasing order is introduced, and a mixed integer quadratic formulation for general weights is presented. In both cases, we show that these formulations allow us the use of kernel functions to construct non linear classifiers.  Besides, we report some computational results about the predictive performance of the introduced approach (also in its kernel version) in comparison with other SVM models existing in the literature. 
\end{abstract}

{\bf Keywords:}
	Data Science; Classification; Support Vector Machine; OWA Operators, Mixed Integer Quadratic Programming.



\section{Introduction}

Support Vector Machine (SVM) models have become one of the most used approaches of Mathematical Programming to address classification problems. SVM techniques have been applied in many different fields since the introduction of the classical soft margin SVM by \cite{13SVN} and \cite{Vapnikunicidad}. Among them, image recognition, bioinformatics, and face detection, see \cite{CERVANTES2020} and references therein.

In the literature, many models based on the classical soft margin SVM approach have been developed with the aim of improving its predictive performance. For instance, different norms have been used to measure the margin between classes, \citep{BlancoPuertoChia}; ouliers or label noise influence has been considered in some models \citep{Baldomero2020, Baldomero2021,BlaJapPue20LN}; and other models also consider feature selection \citep{23Antonio,Jimenez,22Luisa,Lee, 21Maldonado2014}.

In  \cite{MMM18} a new approach for the classical soft-margin Support Vector Machine is proposed. This methodology proposes to apply the OWA operator to modify the hinge loss function of classical SVM. The idea is rather appealing in that it allows to tune the importance of deviations according to their size. Thus, accounting differently classification errors with respect to a preference ranking induced by the sorting of their sizes (distances to the classification hyperplane).

The idea of penalizing the classification errors unevenly according to their sizes is aligned with the use of OWA operators that have become very popular in different areas  of decision theory. Surprisingly, although very natural, this approach has been never tried in SVM and these authors propose a two-step heuristic method to solve their model: 1) the classical SVM is trained and its classification errors induce an order based on the distances to the classifier; 2) the SVM is re-trained using a weighted sum of classification errors with weights induced from the order of the solution in the first step. In \cite{Maldonado2020}, they propose an analogous method, but using
the induced order of fuzzy density-based methods for outlier detection.
This approach is very simple and has the same complexity as the classical SVM beyond of being applied twice. Moreover, as the authors show in their paper, its predictive performance is superior to the traditional SVM in a set of databases that are reported in the paper.

We find the idea of applying OWA operators to classical SVM remarkable and we would like to contribute a bit more on it. Actually, our analysis shows that the approach by \cite{MMM18} is a first step into this direction but it is not really an OWA operator. In fact, we think that there is still room for new contributions and improvements on this topic. Actually, we show in this paper that the methodology by \cite{MMM18}, denoted from now on as app-OWA-SVM, is a heuristic approximation to the exact application of OWA operators to classical SVM. Our aim with this paper is to develop the exact methodology for including OWA operators into SVM and to compare these results with classical SVM and with the approximated OWA version proposed by \cite{MMM18}.

Our contribution in this paper is the following. We present the exact application of OWA operators to the soft-margin hinge loss SVM and prove that the kernel extension of this methodology is also possible. Our analysis  distinguishes between convex OWA operators (those induced by monotone non-decreasing weights) and non-convex ones. For the first family of methods the complexity of the exact OWA-SVM is similar to the classical SVM. However, the second family of methods, namely the non-convex ones, is more complex in that it involves solving mixed-integer second-order cone programs. Yet, kernel transformations are also satisfactorily extended to this class and moreover, the resulting problems can be solved to optimality with MIP-solvers such as Gurobi, Cplex or Xpress. We 
test the performance of OWA-SVM compared with classical SVM and with  app-OWA-SVM. Our results confirm those already reported by \cite{MMM18}: OWA-SVM is superior to SVM and it performs similarly to app-OWA-SVM. This last observation  shows that app-OWA-SVM 
provides good predictive performance similar to the one of the exact OWA-SVM.

The remainder of this paper is structured as follows. In Section \ref{introproblem}, some notation and details about the problem are described. Section \ref{convex} is devoted to the development of an SVM model which includes the OWA when non-decreasing weights are considered. Besides, in Section \ref{nonconvex} we introduce a general mixed integer quadratic model which allows the use of the OWA for general weights (not necessarily non-decreasing). Section \ref{computational} contains computational experiments carried out on several datasets. Finally, Section \ref{conclusions} includes conclusions and some future research lines.

\section{Soft margin hinge loss SVM including OWA operators}\label{introproblem}

In binary classification problems, we are given a  training set of individuals, $N=\{1,\ldots,n\}$, divided into two classes. Each individual, $i$, is represented by a pair $(\bm{ x}_i,y_i)\in\mathbb{R}^d\times \{-1,1\}$, where $d$ is the number of considered features, $\bm{x}_i$ is a vector with features' values and $y_i$ is the label associated with the class of the individual. The goal of SVM models is to determine a hyperplane  $\bm{w}^T \bm{x} + b$  that optimally separates the training set and that allows the classification of new individuals.

The classical soft margin SVM  model is a compromise between maximizing the distance (margin)
between the two parallel class-supporting hyperplanes and minimizing the deviations of misclassified individuals. It is formulated as follows, see \cite{BraMa98},

\begin{alignat}{3}
\mbox{($\ell_2$-SVM)}& \min_{\bm{w},b,\bm{\xi}}&  \quad & \frac{1}{2}\norm{\bm{w}}_2^2+C\sum_{i=1}^n\xi_{i}, & & \nonumber\\
&\mbox{s.t.} &\quad &y_i(\bm{w}^T\bm{x}_i+b)\geq 1-\xi_i,&\quad&i\in N,\nonumber\\
& & & \bm{w}\in \mathbb{R}^d,\nonumber\\
& & &b\in\mathbb{R},\nonumber\\
& & &\xi_i\geq 0,&\quad&i\in N. \nonumber
\end{alignat}  

In this formulation, $\bm{w}$- and $b$-variables are the coefficients of the optimal separator hyperplane and $\bm{\xi}$-variables represent the deviations associated with missclassified individuals. The margin between both supporting hyperplanes is given by $\frac{2}{\norm{\bm{w}}_2}$. Consequently, as mentioned before, the objective function is a balance between the maximization of the margin and the minimization of the deviations. Observe that this balance is regulated by the constant parameter $C$.

Non-linear classifiers can also be obtained by using the classical SVM model. In order to determine a non linear separator, data of the training set $N$ are mapped onto a higher dimension space by using a projection function $\phi(\cdot)$. By the use of duality theory and kernel functions, one can determine the optimal separator without explicitly knowing $\phi(\cdot)$.
To clarify this aspect, it should be mentioned that kernel functions are those such that can be expressed as $K(\bm{x}_i,\bm{x}_j)=\phi(\bm{x}_i)\cdot\phi(\bm{x}_j)$, where $\cdot$ denotes the scalar product. This, together with the fact that 
dual formulation of $\ell_2$-SVM and the resulting optimal separating hyperplane only depend on the dot product of training samples, makes unnecessary the explicit use of $\phi(\cdot)$. For more details about this kernel-based method, see \cite{11burges1998}.

In the context of SVM models, OWA  operators can be applied to the second term of the objective function of the classical SVM, i.e., considering the ordered weighted sum of deviations of misclassified individuals instead of the sum of them. 
The idea of OWA for a set of amounts is to consider the weighted sum of them but taking into account that the weights are assigned depending on the positions in the ordered sequence of these amounts.
For instance, given a deviation vector $\bm{\xi}'$, the ordered weighted sum of the components of this vector is $\sum_{i=1}^n\lambda_i\xi'_{(i)}$. Where $\bm{\xi}'_{()}=(\xi'_{(1)}, \xi'_{(2)}, \ldots,\xi'_{(n)})$ is the vector $\bm{\xi}'$ with its elements sorted in non-decreasing order and $\lambda_i\geq 0$ represents the weight associated with the $i$-th position of the ordered vector $\bm{\xi}'_{()}$. Observe that for the case $\lambda_i=1,$ $\forall i\in N$, we obtain the sum of these amounts.
Hence, a new SVM model considering OWA operator can be expressed as follows,
\begin{alignat}{3}
& \min_{\bm{w},b,\bm{\xi}}&  \quad & \frac{1}{2}\norm{\bm{w}}_2^2+C\sum_{i=1}^n \lambda_i\xi_{(i)}, & & \nonumber\\
&\mbox{s.t.} &\quad &y_i(\bm{w}^T\bm{x}_i+b)\geq 1-\xi_i,&\quad&i\in N,\nonumber\\
& & & \xi_{(i)}\leq \xi_{(i+1)},&&i=1,\ldots,n-1,\nonumber\\
& & & \bm{w}\in \mathbb{R}^d,\nonumber\\
& & &b\in\mathbb{R},\nonumber\\
& & &\xi_i\geq 0,&\quad&i\in N, \nonumber
\end{alignat}

\noindent
where $\xi_{(i)}$ is a variable that represents the $i$-th smallest deviation among the elements in the training set $N$. Note that elements of vector $\bm{\xi}_{()}$ are equal to the ones of $\bm{\xi}$ but sorted in non-decreasing order.

In the next sections, we address the formulation of this problem. Recall that our objective is to provide an exact methodology for dealing with OWA operators with soft margin SVM. For this purpose, we distinguish between convex and non convex OWA operators. The reason of the aforementioned distinction is that, as we will detail in Section \ref{convex}, the use of non-decreasing weights (convex case) allows to build a quadratic continuous formulation whose difficulty is similar to that of $\ell_2$-SVM. In contrast, the use of non-convex OWA operators leads to the introduction of a mixed integer quadratic programming model which is computationally more complex. Section \ref{nonconvex}  deals with the use of these non-convex OWA operators. Besides, in sections \ref{convex} and \ref{nonconvex}, it will also be discussed if non-linear kernels can  be accommodated in each model.  

\section{An SVM-model introducing convex OWA operators}\label{convex}
In order to apply a correct OWA operator to the deviations of the SVM model, one has to multiply sorted deviation by the corresponding $\bm{\lambda}$-weight in the formulation. We begin analyzing the case of monotone non-decreasing $\bm{\lambda}$-weights since, as we will show, it induces simpler mathematical programming models.

With the aim of providing a formulation of this problem, together with the $\bm{w}$-, $b$- and $\bm{\xi}$-variables used in the classical $\ell_2$-SVM, we need to include a new set of variables to model the order of the deviations of misclassified individuals. In particular, we define	
\begin{equation}
z_{ij}=\begin{cases}
1,&\mbox{ if deviation of observation $i$ is in the $j$-th position of the sorted}\\
&\mbox{ vector of deviations,}\\
0,&\mbox{ otherwise,}\\
\end{cases}
\label{zvariables}
\end{equation}
\noindent
for $i,j\in N$. 
Given a vector of the deviation values related to each individual, $\bm{\xi}'$, the use of $\bm{z}$-variables allows us to express the ordered weighted average of these deviations with $\bm{\lambda}$-weights given in non-decreasing order as follows,
\begin{alignat}{3}
& \sum_{i=1}^n\lambda_i\xi_{(i)}'\quad=& \quad\max_{\bm{z}}\quad &  \sum_{i=1}^n\sum_{j=1}^n \lambda_j\xi_i' z_{ij}, & &\nonumber \\
&&\quad\mbox{s.t.}\quad& \sum_{i=1}^n z_{ij}=1,&\quad&j\in N,\label{assign1}\\
&&&\sum_{j=1}^n z_{ij}=1,&\quad&i\in N,\label{assign2}\\
&&& z_{ij}\geq 0,&&i,j\in N.\label{zvar}
\end{alignat} 
Constraints \eqref{assign1} and \eqref{assign2} ensure, respectively, that exactly one element of $N$ is in each position and that each position is allocated to exactly one element of $N$. Besides, due to total unimodularity property, $\bm{z}$-variables can be relaxed as presented in \eqref{zvar}. Hence, a formulation of the SVM with convex OWA operators is given by
\begin{alignat}{3}
& \min_{\bm{w},b,\bm{\xi}}&  \quad & \frac{1}{2}\norm{\bm{w}}_2^2+\max_{\bm{z}} \sum_{i=1}^n\sum_{j=1}^n C\lambda_j\xi_i z_{ij}, & &\nonumber \\
&\mbox{s.t.} &\quad&\eqref{assign1}-\eqref{zvar},\nonumber\\ 
&&&y_i(\bm{w}^T\bm{x}_i+b)\geq 1-\xi_i,&\quad&i\in N, \label{clasconst}\\
& & & \bm{w}\in \mathbb{R}^d,\label{wvar}\\
& & &b\in\mathbb{R},\label{bvar}\\
& & &\xi_i\geq 0,&\quad&i\in N. \label{xivar}
\end{alignat}

Like in the $\ell_2$-SVM formulation, 
constraints \eqref{clasconst} are the classical ones appearing in $\ell_2$-SVM and the restrictions which determine the deviations of misclassified elements of $N$. Constraints \eqref{wvar}-\eqref{xivar} determine the domains of the corresponding variables.

Observe that this optimization model includes an inner maximization problem which intends to obtain the OWA of deviations taking advantage of the fact of using a non-decreasing weight vector. Considering the results of \cite{Blanco14} in the context of facility location problems, we obtain the following quadratic continuous formulation dualizing the inner problem. 

\begin{alignat}{3}
(\convex)\quad& \min_{\bm{w},b,\bm{\xi},\bm{u},\bm{v}}&  \quad & \frac{1}{2}\norm{\bm{w}}_2^2+\sum_{i=1}^n u_i+\sum_{j=1}^n v_j, & & \nonumber\\
&\mbox{s.t.} &\quad &\eqref{clasconst},\eqref{wvar}-\eqref{xivar},&\nonumber\\
&&& u_i+v_j\geq C\lambda_j\xi_i, &&j\in N,\\
&&& u_i \in \mathbb{R},&&i\in N,\\
&&& v_j \in \mathbb{R},&&j\in N.
\end{alignat} 

\noindent Where $\bm{v}$- and $\bm{u}$-variables are dual variables associated with constraints \eqref{assign1} and \eqref{assign2}, respectively. Note that \convex is a quadratic continuous model which determines a linear classifier considering ordered weighted average of individuals errors. 

\begin{remark}
By considering the model proposed in \cite{OGRYCZAK2003} for minimizing the sum of $k$ largest functions, an alternative formulation to \convex\ is

\begin{alignat*}{3}
\mbox{\normalfont{(OT-C-OWA-SVM)}}\quad& \min_{ \bm{w},b,\bm{\xi},\bm{t},\bm{d}}&  \quad & \frac{1}{2}\norm{\bm{w}}_2^2+\sum_{k=1}^n(\lambda_{n-k+1}-\lambda_{n-k})\left(kt_k+\sum_{i=1}^n d_{ik}\right) & &  \\
&\mbox{\normalfont s.t.} &\quad &\eqref{clasconst},\eqref{wvar}-\eqref{xivar},&\nonumber\\
&&&d_{ik}\geq C\xi_i-t_k,&&\hspace{-1.5cm}i,k\in N,\\ 
&&&d_{ik}\geq 0,&&\hspace{-1.5cm}i,k\in N,\\
&&&t_k\in\mathbb{R},&&\hspace{-1.5cm}k\in N.
\end{alignat*}
\noindent Some preliminary computational results show that formulation \convex\ outperforms, in terms of computational times, formulation OT-C-OWA-SVM.
\end{remark}

As in classical SVM, it would be interesting to check whether it is possible to develop a methodology for obtaining non-linear separators by applying the kernel trick. For this reason, once we have a primal formulation of \convex, we present  its dual version that will be very useful to build non linear classifiers. The following results give a formulation of the dual problem.

\begin{prop}
The dual form of \convex is given by:

\begin{alignat}{3}
& \max_{\bm{\alpha},\bm{\eta}}&  \quad & \sum_{i=1}^n\alpha_i-\frac{1}{2}\sum_{i=1}^n\sum_{j=1}^n\alpha_i\alpha_jy_iy_j\bm{x}_i\cdot \bm{x}_j, & &\nonumber \\
&\mbox{\normalfont s.t.} &\quad &\sum_{i=1}^n\alpha_i y_i=0,&&\label{cdual1}\\
&&&\alpha_i\leq \sum_{j=1}^n\eta_{ij}C\lambda_j,&&i\in N,\label{cdual2}\\\
&&&\sum_{i=1}^n\eta_{ij}=1,&&j\in N,\label{cdual3}\\\
&&&\sum_{j=1}^n\eta_{ij}=1,&&i\in N,\label{cdual4}\\\
&&&0\leq \alpha_i,&&i\in N,\label{cdual5}\\\
&&&0\leq\eta_{ij},&&i,j\in N.\label{cdual6}\
\end{alignat} 

\end{prop}

\dem
The Lagrangian function associated with model \convex is

\begin{eqnarray*}
L(\bm{w},b,\bm{\xi},\bm{u},\bm{v})&=&\frac{1}{2}\norm{\bm{w}}_2^2+\sum_{i=1}^n u_i+\sum_{j=1}^n v_j+\sum_{i=1}^n\alpha_i[1-\xi_i-y_i(\bm{w}^T \bm{x}_i+b)]\\
&&-\sum_{i=1}^n\mu_i\xi_i+\sum_{i=1}^n\sum_{j=1}^n \eta_{ij}(C\lambda_j\xi_i-u_i-v_j),
\end{eqnarray*}
where $\bm{\alpha}\geq 0$, $\bm{\mu}\geq 0$ and $\bm{\eta}\geq 0$ are positive Lagrangian multipliers. The necessary and sufficient optimality conditions for \convex result in:

\begin{alignat}{2}
\frac{\partial L (\bm{w},b,\bm{\xi},\bm{u},\bm{v})}{\partial w_j}&=w_j-\sum_{i=1}^n\alpha_iy_ix_{ij}=0,&j\in N,\label{c1}\\
\frac{\partial L (\bm{w},b,\bm{\xi},\bm{u},\bm{v})}{\partial b}&=-\sum_{i=1}^n\alpha_i y_i=0,&\label{c2}\\
\frac{\partial L (\bm{w},b,\bm{\xi},\bm{u},\bm{v})}{\partial \xi_i}&=-\alpha_i-\mu_i+\sum_{j=1}^n\eta_{ij}C\lambda_j=0,&i\in N,\label{c3}\\
\frac{\partial L (\bm{w},b,\bm{\xi},\bm{u},\bm{v})}{\partial u_i}&=1-\sum_{j=1}^n\eta_{ij}=0,&i\in N,\label{c4}\\
\frac{\partial L (\bm{w},b,\bm{\xi},\bm{u},\bm{v})}{\partial v_j}&=1-\sum_{i=1}^n\eta_{ij}=0,&j\in N,\label{c5}\\
\alpha_i[1-\xi_i-y_i(\bm{w}^T\bm{x}_i+b)]&=0,&i\in N,\label{c6}\\
\mu_i\xi_i&=0,&i\in N,\label{c7}\\
\eta_{ij}(C\lambda_j\xi_i-u_i-v_j)&=0,&i,j\in N,\label{c8}\\
\alpha_i,\mu_i,\eta_{ij}&\geq 0,&i,j\in N.\label{c9}
\end{alignat}

From \eqref{c1}, it can be shown that $w_j=\sum_{i=1}^n\alpha_i y_i x_{ij}$ for $j\in N$. Besides, as in the classic SVM, condition \eqref{c2} results in constraint \eqref{cdual1}. In addition, conditions \eqref{c3} can be replaced by inequalities \eqref{cdual2}. Finally, constraints \eqref{cdual3} and \eqref{cdual4} can be deduced from conditions \eqref{c4} and \eqref{c5}, respectively.

By using the complementary slackness conditions and replacing  $w_j=\sum_{i=1}^n\alpha_i y_i x_{ij}$ in the Lagrangian function, we obtain
\begin{eqnarray*}
L(\bm{w},b,\xi,u,v)&=&\sum_{i=1}^n\alpha_i-\frac{1}{2}\sum_{i=1}^n\sum_{j=1}^n\alpha_i\alpha_jy_iy_j\bm{x}_i\cdot\bm{x}_j.
\end{eqnarray*}
\fin
Observe that the formulation of the dual problem depends on the observed data through the scalar product of two observations. 
Hence, replacing in \convex the scalar products of training data by a kernel function $K(\bm{x}_i,\bm{x}_j)=\phi(\bm{x}_i)\cdot\phi(\bm{x}_j)$, where $\phi$ is a map of the observations in a higher dimension space, the resulting formulation is

\begin{alignat}{3}
(\kernelconvex)\quad& \max_{\bm{\alpha},\bm{\eta}}&  \quad & \sum_{i=1}^n\alpha_i-\frac{1}{2}\sum_{i=1}^n\sum_{j=1}^n\alpha_i\alpha_jy_iy_j K(\bm{x}_i,\bm{x}_j), & & \nonumber\\
&\mbox{s.t.} &\quad &\eqref{cdual1}-\eqref{cdual6}.&&\nonumber
\end{alignat} 

Thus, kernel trick can be applied. Recall that this trick consists in using kernel functions in such a way that it is not necessary to explicitly know the transformation $\phi(\cdot)$ and  this formulation does not depend on the dimension of feature space.

Given a sample, $\bm{x}$, belonging to an unknown class, the separator function of a non linear SVM is given by
\begin{equation}
\bm{w}_{\phi}^T \phi(\bm{x})+b=\sum_{i=1}^n\alpha^*_iy_i\phi(\bm{x}_i)\cdot\phi(\bm{x})+b=\sum_{i=1}^n\alpha^*_iy_iK(\bm{x}_i,\bm{x})+b.\label{nonlinear}
\end{equation}
\noindent When using \kernelconvex to obtain a non linear separator, $\bm{\alpha}^*$ values are given by the optimal solution of \kernelconvex. The following result states how to determine the value of $b$ coefficient.

\begin{prop}\label{propb}
$b$-coefficient of the separator function associated with \kernelconvex is given by
$$b=\frac{1-y_k(\sum_{i=1}^n\alpha^*_iy_iK(\bm{x}_i,\bm{x}_k))}{y_k},$$

\noindent	where $k\in N$ verifies that $0<\alpha^*_k<C\sum_{j=1}^n\eta_{kj}^*\lambda_j$, and $(\bm{\alpha}^*,\bm{\eta}^*)$ are the optimal values of \kernelconvex.	
\end{prop}

\dem
Let $k\in N$ be such that  $0<\alpha^*_k<C\sum_{j=1}^n\eta_{kj}^*\lambda_j$. Then, due to conditions \eqref{c3}, $\mu^*_k>0$. Since \eqref{c7} holds, $\xi^*_k=0$. Considering \eqref{c6}, we obtain 
$$1-y_k(\sum_{i=1}^n\alpha^*_i y_i K(\bm{x}_i,\bm{x}_k)+b)=0.$$

Then,
$$b=\frac{1-y_k(\sum_{i=1}^n\alpha^*_iy_iK(\bm{x}_i,\bm{x}_k))}{y_k}.$$
\fin

The discussion above proves that kernel trick can be used in OWA-SVM provided that $\bm{\lambda}$-weights are given in non-decreasing order and that \kernelconvex formulation is used. Furthermore, the non linear separator can be easily obtained by using \eqref{nonlinear} and Proposition \ref{propb}.

Regarding formulation \kernelconvex, we find some differences with respect to the classical kernel extension of $\ell_2$-SVM. Specifically, it is necessary to include some variables ($\bm{\eta}$) and constraints (\eqref{cdual2}-\eqref{cdual4},\eqref{cdual6}) which do not appear in the classical dual model. Despite this, the resulting formulation \kernelconvex is a convex quadratic continuous formulation that can be solved in reasonable small times comparable to the times of classical kernel-version SVM model as we will see in Section \ref{computational}. Then, we have obtained an exact OWA-SVM approach for non-decreasing $\lambda$-weights that can be efficiently solved.

\section{An SVM-model introducing non convex OWA operators}\label{nonconvex}

The main goal of this Section is to introduce OWA in the SVM model when general $\bm{\lambda}$-weights are considered, not necessarily given in non decreasing order. In contrast with the formulation addressed in the previous section, the use of general $\bm{\lambda}$-weights forces the introduction of binary variables in the formulation in order to model the sorting of the SVM related deviations. As a consequence, the resulting model is a quadratic mixed integer formulation which is computationally more complex than \convex. Besides, in spite of using a MIQP to model OWA-SVM with general $\bm{\lambda}$-weights, we  are able to deal with a kernel extension in a different way.

As previously mentioned, OWA operators with general $\bm{\lambda}$-weights have been applied in many combinatorial optimization problems. Particularly, in \cite{FPP14}, OWA problems are analyzed from a modeling point of view and several formulations are compared. Based on this analysis, we present a quadratic mixed integer formulation for the SVM model that we are studying.

To this purpose, it is necessary to use the $\bm{z}$-variables described in \eqref{zvariables} and to introduce a new family of continuous variables, for $k\in N$, defined as:
\begin{eqnarray*}
\theta_{k}&=&\mbox{deviation associated with the individual which is in the $k$-th position}\\
&&\mbox{of the sorted vector of deviations, }
\end{eqnarray*}

The resulting formulation is the following:

\begin{alignat}{3}
(\nonconvex)\quad& \min_{\bm{w},b,\bm{\xi},\bm{z},\bm{\theta}}&  \quad & \frac{1}{2}\norm{\bm{w}}_2^2+C\sum_{k=1}^n\lambda_k\theta_{k}, & & \nonumber\\
&\mbox{s.t.} &\quad &\eqref{clasconst},\eqref{assign1},\eqref{wvar}-\eqref{xivar},\nonumber\\
& & &\theta_{k}\geq\xi_i-M(1-\sum_{\stackrel{j=1}{j\leq k}}^nz_{ij}),&\quad&i,k\in N,\label{the_xi}\\
& & &z_{ik}\in\{0,1\},&&i,k\in N,\label{zbin}\\
& & &\theta_{k}\geq 0,&&k\in N.\label{thetavar}
\end{alignat} 

Constraints \eqref{the_xi} ensure that deviation value in position $k$ is at least the deviation of element $i$, if $i$ is in a position smaller than or equal to $k$, for $i,k\in N$. Constraints \eqref{the_xi} use a big $M$ parameter to establish this link between  $\theta_k$- and $\xi_i$-variables. Note that the maximum distance between two points of the training data is a valid value of $M$.

Obseve that, in contrast with formulation \convex, it is necessary to include binary variables to correctly model the order. As a consequence, completely different techniques must be applied to extend kernel trick to this formulation. 
In what follows, we develop a model to accommodate non linear kernel functions in \nonconvex, see \cite{Brooks2011}. 

\begin{remark}
In model \nonconvex, assume that $\bm{z}$-variables are fixed to $\hat{\bm{z}}$ (feasible assignment). Then the following formulation can be stated,
\begin{alignat}{3}
(\nonconvexz)\quad& \min_{\bm{w},b,\bm{\xi},\bm{\theta}}&  \quad & \frac{1}{2}\norm{\bm{w}}_2^2+C\sum_{k=1}^n\lambda_k\theta_{k}, & & \nonumber\\
&\mbox{\normalfont s.t.} &\quad &\eqref{clasconst},\eqref{wvar}-\eqref{xivar},\eqref{thetavar},\nonumber\\
& & &\theta_{k}\geq\xi_i-M(1-\sum_{\stackrel{j=1}{j\leq k}}^n\hat{z}_{ij}),&\quad&i,k\in N.
\end{alignat}

The dual formulation of \nonconvexz is  
\begin{alignat}{3}
(\nonconvexdual)\quad& \max_{\bm{\alpha},\bm{\mu}}&  \quad & \sum_{i=1}^n\alpha_i-\frac{1}{2}\sum_{i=1}^n\sum_{j=1}^n\alpha_i\alpha_jy_iy_j\bm{x}_i\cdot\bm{x}_j\nonumber\\
&&&-
M\sum_{i=1}^n\sum_{k=1}^n\mu_{ik}\left(1-\sum_{j\leq k}\hat{z}_{ij}\right), & & \nonumber\\
&\mbox{\normalfont s.t.} &\quad &\eqref{cdual1},\eqref{cdual5},\nonumber\\
&&&\alpha_i\leq\sum_{k=1}^n\mu_{ik},&&i\in N,\label{alpha1}\\
&&&\sum_{i=1}^n\mu_{ik}\leq C\lambda_k,&&i,k\in N,\label{alpha2}\\
&&&\mu_{ik}\geq 0,&&\hspace{-4cm}i,k\in N.
\end{alignat}

Besides, from necessary and sufficient optimality conditions, $\displaystyle w_j=\sum_{i=1}^n\alpha_iy_ix_{ij}$ for $j\in N$.

\end{remark}

Based on the link between \nonconvexz and \nonconvexdual,
we propose an alternative formulation of OWA-SVM for general weights. This formulation allows the use of kernel functions, and consequently, the use of non linear separators.

\begin{alignat}{3}
(\nonconvextilde)\quad&\min_{\bm{\alpha},b,\bm{\xi},\bm{\theta},\bm{z}}&  \quad & \frac{1}{2}\sum_{i=1}^n\sum_{j=1}^n y_i y_j\alpha_i\alpha_j\bm{x}_i\cdot\bm{x}_j+C\sum_{k=1}^n\lambda_k\theta_{k}, & & \nonumber\\
&\mbox{s.t.}&&\eqref{assign1},\eqref{bvar},\eqref{xivar},\eqref{the_xi}-\eqref{thetavar},\nonumber\\
&& &y_i(\sum_{j=1}^n y_j\alpha_j\bm{x}_i\cdot\bm{x}_j+b)\geq 1-\xi_i,&\quad&i\in N.
\end{alignat} 

\begin{prop}\label{equivalence1}
Given an optimal solution of \nonconvex, $(\bm{w}^*,b^*,\bm{\xi}^*,\bm{\theta}^*,\bm{z}^*)$, it can be built a feasible solution of \nonconvextilde, $(\bm{\alpha}^*,b^*,\bm{\xi}^*,\bm{\theta}^*,\bm{z}^*)$, with the same objective value.
\end{prop}
\dem

Given a solution of \nonconvex, $(\bm{w}^*,b^*,\bm{\xi}^*,\bm{\theta}^*,\bm{z}^*)$, then $(\bm{w}^*,b^*,\bm{\xi}^*,\bm{\theta}^*)$ is an optimal solution of NC-OWA-SVM($\bm{z}^*$). From necessary and sufficient optimality conditions, $\bm{w}^*=\sum_{i=1}^n\alpha'_i y_i x_{ij}$, where $\bm{\alpha}'$ are the optimal values of $\bm{\alpha}$ variables appearing in NC-OWA-SVM$_D$($\bm{z}^*$). By defining $\bm{\alpha}^*=\bm{\alpha}'$, $(\bm{\alpha}^*,b^*,\bm{\xi}^*,\bm{\theta}^*,\bm{z}^*)$ is feasible for \nonconvextilde and 

$$\frac{1}{2}\norm{\bm{w}^*}_2^2+C\sum_{k=1}^n\lambda_k\theta^*_{k}=\frac{1}{2}\sum_{i=1}^n\sum_{j=1}^n y_i y_j\alpha^*_i\alpha^*_j\bm{x}_i\cdot\bm{x}_j+C\sum_{k=1}^n\lambda_k\theta^*_{k}.$$

\fin
\begin{prop}\label{equivalence2}
Given an optimal solution of \nonconvextilde, $(\bm{\alpha}^*,b^*,\bm{\xi}^*,\bm{\theta}^*,\bm{z}^*)$, a feasible solution of \nonconvex with the same objective value, $(\bm{ w}^*,b^*,\bm{\xi}^*,\bm{\theta}^*,\bm{z}^*)$, can be built.
\end{prop}

\dem

Given an optimal solution of \nonconvextilde, $(\bm{\alpha}^*,b^*,\bm{\xi}^*,\bm{\theta}^*,\bm{z}^*)$, we define 
$w^*_j=\sum_{i=1}^n\alpha^*_i y_i x_{ij}$. This solution is feasible for \nonconvex since \nonconvextilde has been built from \nonconvex, by replacing $w_j$ by $\sum_{i=1}^n\alpha_i y_i x_{ij}$.

\fin

Propositions \ref{equivalence1} and \ref{equivalence2} show that formulations \nonconvex and \nonconvextilde are equivalent in the sense that optimal solutions to \nonconvex can be built from  optimal solutions of \nonconvextilde\  with the same objective value, and vice\-versa.
Note that formulation \nonconvextilde allows us to accommodate non linear kernel functions. Replacing scalar products by a general kernel function, the resulting formulation is

\begin{alignat}{3}
(\kernelnonconvex)\quad &\min_{\bm{\alpha},b,\bm{\xi},\bm{\theta},\bm{z}}&  \quad & \frac{1}{2}\sum_{i=1}^n\sum_{j=1}^n y_i y_j\alpha_i\alpha_jK(\bm{x}_i,\bm{ x}_j)+C\sum_{k=1}^n\lambda_k\theta_{k}, & & \nonumber\\
&\mbox{s.t.}&\quad&\eqref{assign1},\eqref{bvar},\eqref{xivar},\eqref{the_xi}-\eqref{thetavar},\nonumber\\
&& &y_i(\sum_{j=1}^n y_j\alpha_j K(\bm{x}_i,\bm{x}_j)+b)\geq 1-\xi_i,&&\hspace{-0.5cm}i\in N.\label{nlxi}
\end{alignat} 

Observe that, in \kernelnonconvex formulation, a valid value for big $M$ parameter should be determined and, the tighter the formulation, the better the performance. In order to obtain a good estimate for $M$, one can solve the following auxiliary problem:
\begin{alignat}{3}
(\MProblem) \quad& \max_{ \bm{\alpha},b,\bm{\xi},\theta}& \quad & \theta, & & \nonumber\\
&\mbox{s.t.} &\quad &\eqref{bvar},\eqref{xivar},\eqref{cdual1},\eqref{nlxi},\nonumber\\
& & &\theta\geq\xi_i,&\quad&i\in N \label{thetaxi}\\
& & &\frac{1}{2}\sum_{i=1}^n\sum_{j=1}^n y_i y_j\alpha_i\alpha_jK(\bm{x}_i,\bm{x}_j)+C\lambda_n\theta\leq \mbox{UB},&&\label{objbound}\\
& & & 0\leq \alpha_i\leq C\sum_{k}\lambda_k,&&i\in N,\label{alphabounds}\\
& & &\theta\geq 0,&&\label{thetabounds}
\end{alignat} 

\noindent where UB is an upper bound on the optimal value of \kernelnonconvex and $\theta$ is a variable that represents the deviation of the indiviual that is in the $n$-th position of the sorted vector of deviations, the largest one. In \MProblem formulation, constraint \eqref{cdual1} must be satisfied since it appears in NC-OWA-SVM$_D$($\bm{z}'$) for each feasible assignment $\bm{z}'$. In addition, constraints \eqref{nlxi} ensure that $\bm{\alpha}$ solutions satisfy the constraints of \kernelnonconvex\ and constraints \eqref{thetaxi} establish that $\theta$ is the largest deviation. Besides, we include constraint \eqref{objbound} which restrict the objective value of the original problem to be smaller than or equal to a certain upper bound. Finally, the remaining constraints determine the bounds of the problem variables. Note that constraints \eqref{alphabounds} result from the combination of constraints \eqref{alpha1} and \eqref{alpha2} appearing in formulation \nonconvexdual.

In a similar way, upper and lower bounds on $b$-variable could be obtained. Particulary,
\begin{alignat*}{3}
(\bUB)\quad& \max_{ \bm{\alpha},b,\bm{\xi},\theta}& \quad & b, & & \nonumber\\
&\mbox{s.t.} &\quad &\eqref{bvar},\eqref{xivar},\eqref{cdual1},\eqref{nlxi}-\eqref{thetabounds}\nonumber
\end{alignat*}  
provides an upper bound on $b$-variable. Analogously,
\begin{alignat*}{3}
(\bLB)\quad& \min_{ \bm{\alpha},b,\bm{\xi},\theta}& \quad & b, & & \nonumber\\
&\mbox{s.t.} &\quad &\eqref{bvar},\eqref{xivar},\eqref{cdual1},\eqref{nlxi}-\eqref{thetabounds}\nonumber
\end{alignat*} 

\noindent allows us to obtain a lower bound on $b$-variable.

We can conclude that, by using an initial upper bound  on \kernelconvex (UB), and auxiliar problems (\MProblem, \bUB, \bLB),  a valid big M value and bounds on the $b$-variable can be determined. Then, \kernelnonconvex can be solved more efficiently. In Algorithm \ref{kernelnonconvexmethod}, the method for solving \kernelnonconvex is outlined.

\begin{algorithm}[htbp]\label{kernelnonconvexmethod}
\KwData{Training sample composed by a set of $n$ individuals with $d$ features.}
\KwResult{OWA-SVM classifier using non convex weights and a certain kernel function $K(\cdot,\cdot)$.}

Solve the dual form of problem $\ell_2$-SVM with kernel function $K(\cdot,\cdot)$ obtaining a solution $(\bm{\alpha}',b')$. 

Consider the deviations associated with the optimal solution $(\bm{\alpha}',b')$ and sort them in non-decreasing order, obtaining a sorted vector of deviations $\theta'$.

Build a feasible solution for \kernelnonconvex:
$$\mbox{UB}^*:=\frac{1}{2}\sum_{i=1}^n\sum_{j=1}^n y_i y_j\alpha'_i\alpha'_jK(\bm{x}_i,\bm{ x}_j)+C\sum_{k=1}^n\lambda_k\theta'_{k}$$

Solve the problems \bLB and \bUB establishing UB$=$UB$^*$ in constraints \eqref{objbound}. Obtain optimal objective values: $l_b$ from \bLB and $u_b$ from \bUB.

Solve the problem \MProblem establishing UB$=$UB$^*$ in constraints \eqref{objbound} and adding constraint:
\begin{equation}
l_b\leq b\leq u_b.\label{boundsb}
\end{equation}
The optimal solution of \MProblem is denoted as $ub_M$.

Solve \kernelnonconvex including constraints \eqref{boundsb} and using $M=ub_M$. The optimal solutions of \kernelnonconvex are denoted by 
$(\bm{\alpha}^*,b^*,\bm{\xi}^*,\bm{\theta}^*,\bm{z}^*)$.


\caption{Method for solving \kernelnonconvex.}
\end{algorithm}

Next section will be devoted to some computational studies on the different OWA-SVM models.

\section{Computational experiments} \label{computational}

As mentioned in the Introduction, \cite{MMM18} present a heuristic
approach to the use of OWA in the soft-margin SVM, that for the sake of presentation we denote by app-OWA-SVM. The app-OWA-SVM approach is a two-step method. In the first step, the classical soft margin SVM is solved and in the second step, the order induced by this optimal solution on the deviations is used to assign fixed weights to a new soft-margin SVM model.

In this Section, we analyze the performance of our exact approach to OWA-SVM in comparison with app-OWA-SVM and the classical soft-margin SVM. These methods are applied to the six datasets described in Table \ref{tab:datasets}, taken from UCI repository, see \cite{DatasetUCI}. Actually, these datasets were the ones used by \cite{MMM18} to check the validity of their approach, app-OWA-SVM. Observe that Table \ref{tab:datasets} details the complete names of the datasets, the sample size, the number of features and the proportion of each class in the sample. 
\begin{table}[htbp]
\centering
\begin{tabular}{llrrr}
	\toprule
	Label & Complete name  & $n$     & $d$     & Class(\%) \\
	\midrule
	IONO & Ionosphere & 351   & 34    & 64.1/35.9 \\
	WBC   & Wisconsin Breast Cancer & 569   & 30    & 62.7/37.3 \\
	AUS & Australian Credit&690 & 14& 55.5/44.5\\
	DIA & Pima Indians Diabetes&768&8&65.1/34.9 \\
	GC & German Credit&1000&24&70.0/30.0\\
	SPL & Splice &1000&60&51.7/48.3\\
	\bottomrule
\end{tabular}%
\caption{Analyzed datasets from UCI repository}
\label{tab:datasets}%
\end{table}%

Our comparison not only includes the previously mentioned models using linear kernel, but also the influence of Gaussian kernel. Recall that the Gaussian kernel function can be expressed as 
$$K(\bm{x}_i,\bm{x}_j)=exp\left(-\frac{\norm{\bm{x}_i-\bm{x}_j}}{2\sigma^2}\right),$$
where $\sigma>0$ is known as the width parameter, see \cite{Scholkopf}. For that reason, in the reported results, we distinguish between the results of the model with and without the Gaussian kernel.

In order to define the parameters, we follow the same strategy as in \cite{MMM18}.
Specifically, we performed a ten fold cross validation for $C$ and $\sigma$ in $\{2^{-7},2^{-6},\ldots,2^{6}, 2^7\}$. 

Moreover, we analyze four OWA weights based on linguistic quantifiers (see \cite{LUUKKA,YAGER}): basic, quadratic, exponential and trigonometric. For the sake of completeness, we recall the expressions of these quantifiers:
\begin{itemize}
\item Basic quantifier: $Q_b(r)=r^{\widetilde{\alpha}}$, $\widetilde{\alpha}\geq 0$.
\item Quadratic quantifier:  $Q_q(r)=\displaystyle\left(\frac{1}{1-\widetilde{\alpha}(r)^{0.5}}\right)$, $\widetilde{\alpha}\geq 0$.
\item Exponential quantifier: $Q_e(r)=e^{-\widetilde{\alpha}r}$, $\widetilde{\alpha}\geq 0$.
\item Trigonometric quantifier: $Q_r(r)=\arcsin(r\widetilde{\alpha})$, $\widetilde{\alpha}\geq 0$.
\end{itemize}

Considering these quantifiers, the associated weights can be determined by calculating
$$\lambda'_i=Q\left(1-\frac{i-1}{n}\right)-Q\left(1-\frac{i}{n}\right),\mbox{ for }i\in N.$$
Hence, the final weights are given by
$$\lambda_i=\frac{\lambda'_i}{\bar{\lambda'}},$$
where $\bar{\lambda'}$ is the average of $\bm{\lambda}'$-vector.

The specific choice of these weights is motivated since they are the ones reported in \cite{MMM18}. 
Note that the quantifiers related to the weights include a new parameter $\widetilde{\alpha}$ which is also validated in $\widetilde{\alpha}\in\{0.2,0.4,0.6,0.8\}$.

To compare the results of the models, two classification performance metrics are presented: the accuracy (ACC) and the area under the curve (AUC). The accuracy is calculated as 
$$\text{ACC}=\dfrac{TP+TN}{TP+TN+FP+FN},$$ where TP are true positives, TN are true negatives, FP false positives and FN false negatives. The area under the curve is given by $$\text{AUC}=\dfrac{\dfrac{TP}{TP+FN}+\dfrac{TN}{TN+FP}}{2}.$$

Regarding the solution methods used for solving the models,
classical soft-margin $\ell_2$-SVM and the $\ell_2$-SVM model using Gaussian kernel are solved with {\sl SVC} function of Scikit Learn module in Python, see \cite{scikit-learn}. Moreover, 
app-OWA-SVM and app-OWA-SVM$_K$, the models with linear and Gaussian kernel (respectively) appearing in \cite{MMM18} are solved by using the two-step method proposed by them.

The exact OWA-SVM models that we propose can be solved depending on the weights with different approaches. Specifically, \kernelconvex is used  for the weights based on basic and exponential quantifiers since they are monotone non-decreasing. Besides, the weight based on the quadratic quantifier is also monotone non-decreasing for $\widetilde{\alpha}=0.2$. For this reason, \kernelconvex is also used with these weights. Note that \kernelconvex\   formulation is also the one used in the linear kernel case, i.e., $K(\bm{x}_i,\bm{x}_j)=\bm{x}_i\cdot \bm{x}_j$.
For the remaining weights, \nonconvex and \kernelnonconvex (following Algorithm \ref{kernelnonconvexmethod}) are applied to obtain the classifiers. 
It should be noted that all our computational studies were performed using CPLEX 20.1.0 in Python on an Intel(R)
Xeon(R) W-2245 CPU 256 GB RAM computer.

Before comparing the predictive performance of the models, we would like to emphasize the level of adequacy provided by app-OWA-SVM with respect to the correct final ranking of the deviatons. Table \ref{tab:order} reports an illustrative example of how the order of deviations behave in app-OWA-SVM methods compared with the correct final ranking. Specifically, Table \ref{tab:order} reports, for the considered datasets, the average percentage of coincidences between the positions that occupy the individuals in the sorted deviations vector of step one and their positions in step two of the app-OWA-SVM method. In addition, in the third column of Table \ref{tab:order}, we report the average percentage of coincidences between the final order in the app-OWA-SVM method and the order provided by \kernelconvex. It should be highlighted that these results are obtained when applying the ten fold cross validation to the models with $C=1$, $\sigma=1$, basic quantifier weight, and $\widetilde{\alpha}=0.6$.

Results of Table \ref{tab:order} show how the induced order of classical SVM, in step 1 in \cite{MMM18}, is not the same as the one resulting in the second step of app-OWA-SVM. This order is neither the same as the sorting obtained by applying \kernelconvex. In contrast with the approaches in \cite{MMM18}, the exact OWA-SVM models, proposed in this paper, set the order of the deviations while solving the model itself and therefore they actually apply exact OWA operators to SVM. This shows that the method in \cite{MMM18} is a OWA-like approach but actually, it is not an exact application of OWA operators.

\begin{table}[htbp]
\centering
\begin{tabular}{lrr}
	\toprule
	Data  & \multicolumn{1}{l}{Step 1 - Step 2 (\%)} & \multicolumn{1}{l}{Step 2 - \kernelconvex (\%)} \\
	\midrule
	IONO  & 1.68\% & 2.06\% \\
	WBC   & 1.17\% & 0.18\% \\
	AUS   & 0.86\% & 1.71\% \\
	DIA   & 0.56\% & 1.62\% \\
	GC    & 0.40\% & 0.36\% \\
	SPL   & 0.03\% & 0.00\% \\
	\bottomrule
\end{tabular}%
\caption{Average percentage of coincidence in the positions of the sorted deviations vector}
\label{tab:order}
\end{table}%

Focusing on the performance classification metrics, Table \ref{tab:acc} reports the best results in terms of ACC provided by the different models. Note that $\ell_2$-SVM and $\ell_2$-SVM$_K$ show the ACC of the classic model with linear and Gaussian kernels, respectively. The column corresponding to app-OWA-SVM presents the results of the model proposed in \cite{MMM18} and column app-OWA-SVM$_K$ shows the results for the model in \cite{MMM18} using the Gaussian kernel. Finally, ex-OWA-SVM and ex-OWA-SVM$_K$ report the best results of our proposed exact OWA-SVM methods using the formulations  in sections \ref{convex} and \ref{nonconvex}. Observe that, for all datasets, the best ACC are either the one provided by the method presented in \cite{MMM18} using the Gaussian kernel, or the ACC of the exact OWA-SVM  model using the Gaussian kernel.

Particularly, ex-OWA-SVM$_K$ provides the best ACC results for WBC and AUS datasets; app-OWA-SVM$_K$ and ex-OWA-SVM$_K$ seem to yield the same results in the IONO and GC cases; whereas the best results for DIA and SPL are obtained by app-OWA-SVM$_K$. In general, we can observe that the accuracy (ACC) of app-OWA-SVM$_K$ and ex-OWA-SVM$_K$ are similar. This indicates that both approaches are worthy in the sense that they improve this classification performance metric with respect to the classical SVM achieving almost the same value.

\begin{table}[htbp]
\renewcommand{\arraystretch}{1.2}
\begin{adjustbox}{max width=1.00\textwidth}
	\centering
	\begin{tabular}{lrrrrrr}
		\toprule
		\multicolumn{7}{c}{ACC (\%)} \\
		\midrule
		Data  & \multicolumn{1}{l}{$\ell_2$-SVM} & \multicolumn{1}{l}{app-OWA-SVM} & \multicolumn{1}{l}{ex-OWA-SVM} & \multicolumn{1}{l}{$\ell_2$-SVM$_K$} & \multicolumn{1}{l}{app-OWA-SVM$_K$} & \multicolumn{1}{l}{ex-OWA-SVM$_K$} \\
		\midrule
		IONO  & 90.60\% & 90.89\% & 90.89\% & 95.44\% & \textbf{95.72\%} & \textbf{95.72\%} \\
		WBC & 98.07\% & 98.07\% & 97.71\% & 98.24\% & 98.42\% & \textbf{98.77\%} \\
		AUS   & 85.51\% & 86.09\% & 85.51\% & 86.38\% & 87.25\% & \textbf{87.39\%} \\
		DIA   & 77.60\% & 77.73\% & 77.73\% & 77.34\% & \textbf{78.38\%} & 78.12\% \\
		GC    & 76.90\% & 77.50\% & 77.30\% & 77.30\% & \textbf{77.50\%} & \textbf{77.50\%} \\
		SPL   & 81.30\% & 82.00\% & 81.70\% & 88.40\% & \textbf{89.80\%} & 89.40\% \\
		\bottomrule
	\end{tabular}%
\end{adjustbox}
\caption{Best ACC results of each model}
\label{tab:acc}%
\end{table}%

Table \ref{tab:auc} reports the AUC for the best combination of parameter values in each case. As for the ACC measure, the best AUC results are always provided by app-OWA-SVM$_K$ and ex-OWA-SVM$_K$. Regarding the results, 
we observe that the AUC of both approaches are quite similar. For the IONO, DIA and SPL datasets, the same AUC is achieved by app-OWA-SVM$_K$ and ex-OWA-SVM$_K$; ex-OWA-SVM$_K$ provides the best results for the WBC and AUS datasets; and finally, the app-OWA-SVM$_K$ reports the best values of AUC for GC dataset. Both approaches, app-OWA-SVM$_K$ and ex-OWA-SVM$_K$, improve the AUC of the classical $\ell_2$-SVM and $\ell_2$-SVM$_K$.

\begin{table}[htbp]
\renewcommand{\arraystretch}{1.2}
\begin{adjustbox}{max width=1.00\textwidth}
	\centering
	\begin{tabular}{lrrrrrr}
		\toprule
		\multicolumn{7}{c}{AUC (\%)} \\
		\midrule
		Data  & \multicolumn{1}{l}{$\ell_2$-SVM} & \multicolumn{1}{l}{app-OWA-SVM} & \multicolumn{1}{l}{ex-OWA-SVM} & \multicolumn{1}{l}{$\ell_2$-SVM$_K$} & \multicolumn{1}{l}{app-OWA-SVM$_K$} & \multicolumn{1}{l}{ex-OWA-SVM$_K$} \\
		\midrule
		IONO  & 88.67\% & 88.89\% & 88.89\% & 94.55\% & \textbf{95.15\%} & \textbf{95.15\%} \\
		WBC & 97.70\% & 97.70\% & 97.14\% & 97.84\% & 98.08\% & \textbf{98.44\%} \\
		AUS   & 86.20\% & 86.67\% & 86.20\% & 86.54\% & 87.46\% & \textbf{87.71\%} \\
		DIA   & 72.32\% & 72.96\% & 73.00\% & 72.19\% & \textbf{73.30\%} & \textbf{73.30\%} \\
		GC    & 68.98\% & 71.00\% & 69.19\% & 68.74\% & \textbf{71.67\%} & 69.90\% \\
		SPL   & 81.41\% & 82.02\% & 81.78\% & 88.41\% & \textbf{89.86\%} & \textbf{89.86\%} \\
		\bottomrule
	\end{tabular}%
\end{adjustbox}
\caption{Best AUC results for each model}
\label{tab:auc}%
\end{table}%

To conclude the analysis, we wish to include some information on the CPU times needed to solve these problems.
Table \ref{tab:times} reports the average solving time (in seconds) per fold of the models for the parameter values that provide the best AUC. For the approach in \cite{MMM18}, we report the time required by the two steps that are involved in the method. 

Concerning the times of ex-OWA-SVM and ex-OWA-SVM$_K$, we note in passing that they correspond to formulation \kernelconvex since, in all cases tested, the best results in terms of accuracy and AUC are obtained using monotone non-decreasing weights. 
Table \ref{tab:times} shows that the exact approach for OWA-SVM requires more time than the classical SVM and also more than the methods presented in \cite{MMM18}. This is due to the fact that models in \cite{MMM18} have essentially the same complexity as the classical SVM.
\begin{table}[htbp]
\renewcommand{\arraystretch}{1.2}
\begin{adjustbox}{max width=1.00\textwidth}
	\centering
	\begin{tabular}{lrrrrrrrr}
		\toprule
		\multicolumn{9}{c}{Time (s)} \\
		\midrule
		\multicolumn{1}{c}{\multirow{2}[4]{*}{Data}} & \multicolumn{1}{c}{\multirow{2}[4]{*}{$\ell_2$-SVM}} & \multicolumn{2}{c}{app-OWA-SVM} & \multicolumn{1}{c}{\multirow{2}[4]{*}{ex-OWA-SVM}} & \multicolumn{1}{c}{\multirow{2}[4]{*}{$\ell_2$-SVM$_K$}} & \multicolumn{2}{c}{app-OWA-SVM$_K$} & \multicolumn{1}{c}{\multirow{2}[4]{*}{ex-OWA-SVM$_K$}} \\
		\cmidrule{3-4}\cmidrule{7-8}          &       & \multicolumn{1}{c}{Step 1} & \multicolumn{1}{c}{Step 2} &       &       & \multicolumn{1}{c}{Step 1} & \multicolumn{1}{c}{Step 2} &  \\
		\midrule
		IONO  & 0.023 & 0.014 & 0.212 & 4.725 & 0.006 & 0.003 & 1.900 & 4.967 \\
		WBC & 0.003 & 0.003 & 0.494 & 12.628 & 0.003 & 0.002 & 4.915 & 9.345 \\
		AUS   & 0.011 & 0.425 & 0.297 & 19.900 & 0.014 & 0.014 & 8.105 & 26.609 \\
		DIA   & 0.013 & 0.009 & 0.286 & 25.133 & 0.013 & 0.008 & 9.944 & 21.536 \\
		GC    & 0.033 & 0.016 & 0.560 & 41.788 & 0.024 & 0.025 & 17.090 & 31.707 \\
		SPL   & 0.022 & 0.522 & 0.815 & 40.561 & 0.057 & 0.038 & 16.387 & 30.866 \\
		\bottomrule
	\end{tabular}%
\end{adjustbox}
\caption{Times per fold of each model for the best parameter values}
\label{tab:times}%
\end{table}%

The aforementioned results show that the exact OWA-SVM models introduced in this paper improve the classification performance metrics of the classical SVM. Furthermore, these measures are similar to the ones reported in previous approaches to OWA-SVM models in terms of ACC and AUC, although they have the advantage of actually capturing the essence of OWA in its application to SVM.

\section{Conclusions}\label{conclusions}

OWA operators have been applied to different problems of decision theory. This paper proposes an exact approach that allows the introduction of OWA operators for the deviation errors appearing in the soft-margin SVM. For this aim, we have distinguished between OWA-SVM with non-decreasing $\bm{\lambda}$-weights and OWA-SVM with general $\bm{\lambda}$-weights (not necessarily non-decreasing). 

The use of non-decreasing weights allowed to formulate the model using only continuous variables. As consequence, a quadratic continuous formulation was developed, \convex. In addition, it was shown that non linear kernels could be accommodated in this model by the use of the dual formulation. The required solving time of this formulation is similar to the classical one.	
In contrast, the use of non-monotone weights in the OWA-SVM leads us to the use of binary variables to model the order of the deviations vector. Then, a mixed integer quadratic formulation, \nonconvex, is necessary to solve this problem. Hence, \nonconvex is more complex than \convex. Despite this, it is also possible to apply non-linear kernel by the use of an alternative formulation.

We have compared the proposed models with the heuristic approach to OWA-SVM appearing in \cite{MMM18}. Regarding the results, we can conclude that app-OWA-SVM provides solutions very different from the optimal solutions of the resulting model of applying actual OWA operators to SVM. However, both methodologies show similar predictive performance improving the ones obtained with the classical SVM (linear and non linear).

As future research, it could be interesting to analyze the OWA-SVM model using other $\ell_p$-norms. Besides, it could also be studied the integration of other aspects to OWA-SVM models such as feature selection, outlier detection or presence of label noise.

\section*{Acknowledgements}
Alfredo Mar\'in has been supported by Spanish Ministry of Science and Innovation under project PID2019-110886RB-I00. Part of this research was conducted while he was on sabbatical at Universidad de Sevilla, Spain.
Luisa I. Martínez-Merino, Antonio M. Rodríguez-Chía and Justo Puerto thank the Agencia Estatal de Investigación (AEI) and the European Regional Development’s funds (ERDF): project MTM2016-74983-C2; 
Regional Government of Andalusia: projects FEDER-UCA18-106895, FEDER-US-1256951, and P18-FR-1422; and Fundación BBVA: project NetmeetData (Ayudas
Fundación BBVA a equipos de investigación científica 2019). Luisa I. Martínez-Merino has been supported by Subprograma Juan de la Cierva Formación 2019 (FJC2019-039023-I).


\end{document}